\documentclass[12pt]{amsart} 
 
\usepackage{amssymb}
\usepackage{epsf}
\usepackage{amsmath,epsf}
\usepackage[latin1]{inputenc}
\usepackage{amsfonts}

\numberwithin{equation}{section}
\newtheorem{thm}{Theorem}[section]
\newtheorem{prop}[thm]{Proposition}
\newtheorem{conj}[thm]{Conjecture}
\newtheorem{lem}[thm]{Lemma}

\newtheorem{rem}[thm]{Remark}

\newenvironment{Remark}{\begin{rem}\sl}{\end{rem}}

\newcommand{\twolines}[2]{{\scriptstyle #1 \atop \scriptstyle
\vphantom{\sum}#2}}

\def\C{{\mathbb C}}
\def\N{{\mathbb N}}

\def\Q{{\mathbb Q}}

\def\S{{\mathcal S}}
\def\RS#1{{\mathbb S}_{#1}}
\def\Young#1{\vbox{\smallskip\offinterlineskip
    \halign{&\vbox{##}\kern-\Thickness\cr #1}}}

\newdimen\Squaresize \Squaresize=20pt
\newdimen\Thickness \Thickness=.1pt
\newdimen\Correction \Correction=7pt

\def\Vide#1{\hbox{
       \vbox to \Squaresize{\vss
          \hbox to \Squaresize{\hss#1 \hss}\vss}
    \hskip-\Correction}
   \kern-\Thickness}

\def\Box#1{\Carre{$#1$}}

\def\Carre#1{\hbox{\vrule width \Thickness
   \vbox to \Squaresize{\hrule height \Thickness\vss
      \hbox to \Squaresize{\hss$\scriptstyle#1$\hss}
   \vss\hrule height\Thickness}
   \unskip\vrule width \Thickness}
   \kern-\Thickness}

\title[Lattice diagram polynomials]{Lattice Diagram polynomials\\  in one
set of variables}

\author[J.-C.~Aval, F.~Bergeron and N.~Bergeron]
{J.-C.~Aval, F.~Bergeron and N.~Bergeron}

\address[Jean-Christophe Aval]
{Laboratoire A2X\\ Universit\'e Bordeaux 1\\ 351 cours de la
Lib\'eration\\ 33405 Talence cedex\\ FRANCE}

\address[Fran\c cois Bergeron]
{D\'epartement de Math\'ematiques\\ Universit\'e
du Qu\'ebec \`a Montr\'eal\\ Montr\'eal, Qu\'ebec, H3C 3P8, CANADA.}

\address[Nantel Bergeron]
{Department of Mathematics and Statistics\\ York University\\
   To\-ron\-to, Ontario M3J 1P3\\ CANADA}

\email[Jean-Christophe Aval]{aval@math.u-bordeaux.fr}
\email[Fran\c cois Bergeron]{bergeron.francois@uqam.ca}
\email[Nantel Bergeron]{bergeron@mathstat.yorku.ca}
\urladdr[Fran\c cois Bergeron]{http://bergeron.math.uqam.ca}
\urladdr[Nantel Bergeron]{http://www.math.yorku.ca/bergeron}

\date{\today}
\thanks{F. Bergeron and N. Bergeron are supported in part by NSERC and FCAR}
\thanks{J.-C. Aval is partially supported by the Conseil R\'egional
d'Aquitaine}


\begin{document}

\begin{abstract}
The space $M_{\mu/i,j}$ spanned by all partial derivatives of the lattice
polynomial
$\Delta_{\mu/i,j}(X;Y)$ is investigated in~\cite{Berg et al} and many
conjectures are given. Here, we prove all these conjectures for the $Y$-free component $M_{\mu/i,j}^0$ of
$M_{\mu/i,j}$.
In particular,  we give an explicit bases for
$M_{\mu/i,j}^0$ which allow us to prove directly the central {\sl four term recurrence}
for these spaces. 
\end{abstract}

\maketitle

\centerline{\bf This paper is dedicated to the memory of Rodica Simion}
\section{Introduction}
We are going to explicitly describe certain $\S_n$-modules of polynomials, in $n$ variables $x_1,\ldots ,x_n$, that are closely
related to classical harmonic polynomials for the symmetric group $\S_n$. These
last polynomials can be characterized by the fact that they satisfy the conditions
   \begin{equation}\label{harmon}
      \sum_{i=1}^n \partial_{x_i}^k P(x_1,\ldots,x_n) =0,\qquad k=1,2,3,\ldots\,.
   \end{equation}
A classical result of Steinberg states that
the set $M_n$, of all harmonic polynomials for the symmetric group, is  
   $$M_n := {\mathcal L}_\partial[\Delta_n],$$
where $\Delta_n=\prod_{i<j} (x_i-x_j)$ is the Vandermonde determinant, and ${\mathcal L}_\partial[\Delta_n]$ denotes
the linear span of all partial derivatives of $\Delta_n$. It is certainly striking to notice that the dimension of $M_n$ is
$n!$, and there is a lot of other nice results related to $M_n$ and its generalization to reflection groups. The spaces studied
here are natural generalizations of these spaces and spaces studied by DeConcini and Procesi in \cite{deconcini} and Garsia and
Procesi
\cite{garsia procesi}.

The point of departure of this
work consists in replacing $\Delta_n$ by natural generalizations of the Vandermonde determinant. 
To this end, let us define a general {\sl lattice diagram} to be any finite subset of $\N\times\N$. The case corresponding
to diagrams of partitions is of special interest. Recall that a {\sl partition} 
$\mu$ of $n$ is a sequence  $\mu_1\geq \mu_2\geq \cdots \geq \mu_k>0$ of decreasing positive
integers such that $n=\mu_1+\cdots +\mu_k$. We denote $\mu\vdash n$ this fact. For  $\mu\vdash n$, we set 
   $$\mu ! := \mu_1!\, \mu_2!\, \cdots \mu_k!\,.$$
The lattice diagram associated to 
a partition $\mu$ is defined to be the set 
  $$\{(i,j)\,:\,0\le i\le k-1,\,
0\le j\le \mu_{i+1}\},$$
 and we use the symbol
$\mu$  both for the partition and its diagram.
Most definitions and conventions used in this text are those of \cite{Berg et al}.
For example, the diagram of the partition $(4,2,1)$  is geometrically represented as
  $$\Young{\Carre{2,0}\cr
           \Carre{1,0}&\Carre{1,1}\cr
           \Carre{0,0}&\Carre{0,1}&\Carre{0,2}&\Carre{0,3}\cr
           }$$
and it consists of the lattice {\sl cells\/}
 $$\{(0,0),(0,1),(0,2),(0,3),(1,0),(1,1),(2,0)\}.$$
Thus, the {\sl coordinates\/} $(r,c)$ of a cell are such that $r+1$ (resp. $c+1$) is the {\sl row number\/} (resp. {\sl column
number\/}) for the cell position, counting  from bottom up (resp. left to right).

Given a lattice diagram $D=\{(r_1,c_1), (r_2,c_2),\ldots , (r_n,c_n)\}$
we define the {\sl lattice determinant}
  $$
  \Delta_D(X;Y):= \det \big\| x_i^{r_j}y_i^{c_j}\big\|_{i,j=1}^n\,,
  $$
where $X=x_1,x_2,\ldots,x_n$ and $Y=y_1,y_2,\ldots,y_n$.
The polynomial
$\Delta_D(X;Y)$ is  bihomogeneous of degree $|r|=r_1+\cdots +r_n$ in $X$
and  degree $|c|=c_1+\cdots +c_n$  in  $Y$. To insure that this definition
associates a unique polynomial to $D$, we order lattice cells in
increasing lexicographic order.

We will need a few more definitions regarding partitions and diagrams. 
For an $n$-cell diagram $D$, a {\sl tableau} of shape $D$ is an injective map $T:D\to
\{1,2,\ldots,n\}$.  If $T(r,c)=m$, we say that $h_T(m):=r$ is the {\sl height\/} of $m$ in $T$.
We say that $T$ is {\sl row increasing} if $T(i,j)<T(k,j)$ whenever
 $i<k$ (when this has a meaning). Similarly we define {\sl column increasing} tableaux, and a {\sl standard}
tableau is one that is both row and column increasing. We denote by
$\RS{D}$ the set of all standard tableau of
shape $D$. 
 
With these definitions out of the way, we come to our object of interest, namely the modules
 $$M_D := {\mathcal L}_\partial[\Delta_D],$$ 
where $D$ is some lattice diagram. The case where $D=1^n$ ($1^n$ is the partition of $n$ with all parts equal to 1) corresponds to
the classical module of harmonic polynomials. This generalization was first considered by Garsia and Haiman~\cite{Garsia Haiman},
in the special case when
$D=\mu$ is the lattice diagram of a partition. Since then, several other cases have been studied (see \cite{Berg et
al} and \cite{Berg2}).  For any
$n$-cell lattice diagram
$D$, the space $M_D$ affords the structure of an $\S_n$-module, through the action of the  symmetric group on polynomials
consisting in permuting variables.  More precisely, a permutation $\sigma\in \S_n$ acts {\sl diagonally\/} on a polynomial
$P(X;Y)$, as follows
  $$\sigma  P(X;Y)\,=\, P(x_{\sigma_1},x_{\sigma_2},\ldots
,x_{\sigma_n};y_{\sigma_1},y_{\sigma_2},\ldots ,y_{\sigma_n}).$$
Under this action, $\Delta_D=\Delta_D(X;Y)$ is clearly an alternant and the action commutes with partial derivatives,
hence
$M_D$ is an invariant subspace of
$\Q[X,Y]$.

Moreover, since $\Delta_D$ is  bihomogeneous, $M_D$
affords the following natural
bigrading. Denoting by ${\mathcal H}_{r,s}[M_D]$ the subspace consisting of the
bihomogeneous elements of degree  $r$ in $X$ and degree $s$ in $Y$, we have
the direct sum
decomposition
  $$
  M_D = \bigoplus_{r=0}^{|p|}  \bigoplus_{s=0}^{|q|} {\mathcal H}_{r,s}[M_D].
  $$
The {\sl bigraded  Hilbert series} of $M_D$ is
  $$
  F_D(q,t) = \sum_{r=0}^{|p|}  \sum_{s=0}^{|q|}  t^r  q^s  \dim
  {\mathcal H}_{r,s}[M_D],
  $$
and the
bigraded character of
$M_D$ is encoded by the following symmetric function,
  $$
  H_D(X;q,t)= \sum_{r=0}^{|p|}  \sum_{s=0}^{|q|}  t^r  q^s \, {\mathcal
F}\big(\hbox{ch}
  {\mathcal H}_{r,s}[M_D]\big), 
  $$
where $\hbox{ch}{\mathcal H}_{r,s}[M_D]$ denotes the character of
${\mathcal H}_{r,s}[M_D]$
and ${\mathcal F}$ is the Frobenius correspondence which maps the
irreducible character
$\chi^\lambda$  to the Schur function $s_\lambda$. We say that this symmetric function is
the {\sl Bigraded Frobenius Characteristic\/} of $M_D$.
 
The $n!$-factorial conjecture of Garsia and Haiman states that for any partition
diagram $\mu$, $H_\mu(X;q,t)$ is none other than a renormalized version of the
Macdonald polynomial associated to $\mu$. This has recently been shown by Haiman~\cite{haiman} using an algebraic geometry
approach. It develops that a very natural and combinatorial recursive approach to the
$n!$-conjecture involves diagrams obtained by removing a single
cell from a partition diagram. It is conjectured in
\cite{Berg et al} that, for such diagrams, the space 
$M_D$ is a direct sum of left regular representations of $\S_n$. More precisely, if
$\mu$ is a partition of
$n+1$, we  denote by
$\mu/ij$ the lattice diagram obtained by
removing one of its cell $(i,j)$ from the diagram $\mu$. We  refer to the cell
$(i,j)$ as
the {\sl hole} of  $\mu/ij$. The conjecture in question states that the number of copies of
the left regular
representations in $M_{\mu/ij}$ is equal to the cardinality of the
$(i,j)$-{\sl shadow}, that is
the cardinality of
$\{(s,t)\in\mu : s\ge i,\,t\ge j\}$. 

This, and more, is all encoded in the
following four term recurrence for the bigraded Frobenius characteristic
$H_{\mu/ij}$ of $M_{\mu/ij}$.

\begin{conj}[\cite{Berg et al}] \label{conj:rec}
For all $(i,j)\in \mu$, we have  $H_{\mu/ij}=C_{\mu/ij}$.
\end{conj}

\noindent Where, $C_{\mu/ij}$ is defined by the following the following ``{\em four term}'' recurrence
   \begin{equation}\label{quatres termes}
   C_{\mu/ij}= {t^\ell-q^{a+1}\over t^\ell-q^a}   C_{\mu/i,j+1}  +
               {t^{\ell+1}-q^{a}\over t^\ell-q^a}  C_{\mu/i+1,j} -
               {t^{\ell+1}-q^{a+1}\over t^\ell-q^a}  C_{\mu/i+1,j+1},
   \end{equation}
where $\ell$  and $a$  give the number of cells that are respectively north
and east of
$(i,j)$ in  $\mu$. As boundary conditions, we set
$C_{\mu/i,j+1}$, $C_{\mu/i,j+1}$  or $C_{\mu/i,j+1}$
 equal to zero when the corresponding
cells $(i,j+1),$ $(i+1,j)$ or  $(i+1,j+1)$
fall outside of $\mu$. Furthermore, if $(i,j)$ is a corner of $\mu$, then $\mu/ij$
is a partition diagram $\nu$,
and we set 
   $$C_{\mu/ij}=H_{\nu}.$$

For any lattice diagram $D$, we consider $M_D^0$  the  $Y$-free component of $M_D$, this is to say that
  $$
  M_D^0 = \bigoplus_{r=0}^{|p|}  {\mathcal H}_{r,0}[M_D]\,.
  $$
In this paper, we study the spaces $M_{\mu/ij}^0$, and show that the $Y$-free specialization  of the
conjecture~\ref{conj:rec} hold for these spaces. In particular this implies that
   \begin{equation}\label{dimension}
      \dim M_{\mu/i,j}^0 ={\frac {n!} {\mu!}} \ |\{\, (r,c)\in\mu\ |\, i\leq r\leq \ell  \},
   \end{equation}
where $\ell$ is the largest integer for which the corresponding row of $\mu$ has at least $j$ cells. Moreover,
we obtain a formula for $H_{\mu/ij}^0$, the graded Frobenius character of $M_{\mu/ij}^0$.

\section{A Basis for $M_\mu^0$} 

In preparation for our description of the modules $M_{\mu/ij}^0$, we need to recall and reformulate some results about the
modules $M_\mu^0$. Although a recursive description of a basis for
$M_\mu^0$ is given in~\cite{nantel adriano}, and a direct description is given in~\cite{aval}, we give here a new description
directly in term of standard tableaux of shape $\mu$. One can immediately link this description to Tanisaki's construction~\cite{Tani} of the
defining ideal of
$M_\mu^0$. Moreover, we will see that it clearly generalizes the ``Artin'' basis for $M_n=M_{1^n}^0$:
  $${\mathcal B}_n:=\{\partial_X^{\bf a}\,\Delta_n(X)\ |\ {\bf a}=(a_1,a_2,\,\ldots\,,a_n), \quad a_i< i\ \}\,,$$
where we use the vectorial notation 
$\partial_ X^{\bf a}:=\partial_{x_1}^{a_1}\,\partial_{x_2}^{a_2}\,\ldots\,\partial_{x_n}^{a_n}$.

First, let $T$ be a tableau of shape $D$ (any diagram), and define
   $$Y_T:=\prod_{(r,c)\in D} y_{_{T(r,c)}}^c,$$
so that, for $1\leq k\leq n$, the exponent of $y_k$ is $c$, if $(r,c)$ is the position of $k$ in $T$. For example,
if 
   $$T=\lower15pt\hbox{\Young{\Box{3}&\Box{6}\cr\Box{1}&\Box{2}&\Box{4}&\Box{5}&\Box{7}\cr}}$$
then $Y_T=y_2\,y_6\,y_4^2\,y_5^3\,y_7^4$. Clearly, the monomial $Y_T$ encodes the columns in which the entries of $T$ appear.
More precisely, let the {\sl column set\/} of a diagram $T$ of shape $D$, be
   $$\Gamma(T):=\{(k,\Gamma_k(T)\ |\ 0\leq k\leq \max_{(r,c)\in D} c\},$$
where for a given $k$, the set $\Gamma_k(T)$ is the set of entries in column $k$, that is
   $$\Gamma_k(T):=\{ \ T(j,k)\ |\ (j,k)\in D\ \}.$$
Then $Y_T=Y_R$, if and only if $T$ and $R$ have same column sets (up to the addition of pairs of the form $(k,\{\})$). Observe that
$T$ and
$R$ need not be of the same shape.

We will denote
${\partial Y}_T$ the differential operator obtained by replacing each $y_k$ in $Y_T$ by $\partial_{y_k}$. With the same
conventions, define
   $$X_T:=\prod_{(r,c)\in D} x_{_{T(r,c)}}^r\,,\qquad {\rm and}\qquad Z_T:=X_T\,Y_T.$$
It is easy to see that $T\mapsto Z_T$ establishes is a bijection between injective tableaux of shape $D$ and monomials of the
form
 $$x_{\sigma_1}^{r_1}\,y_{\sigma_1}^{c_1} \,x_{\sigma_2}^{r_2}y_{\sigma_2}^{c_2}\,\cdots\,
         x_{\sigma_n}^{r_n}y_{\sigma_n}^{c_n},$$
with $\sigma$ in $\S_n$ and $D=\{(r_1,c_1), (r_2,c_2),\ldots , (r_n,c_n)\}$ ordered lexicographically.
If we set
   $$\gamma_D:=\prod_{(r,c)\in D} c!,$$
then we easily get the following lemma.

\begin{lem}\label{lem:tabl} For two tableaux $T$ and $R$, both of shape $D$, one has
  $${\partial Y}_T Z_R=\gamma_D\ X_R$$
if $R$ can be obtained from $T$ by a column fixing permutation of its entries. 
Otherwise 
  $${\partial Y}_T Z_R=0.$$
\end{lem}

Given a tableau
$T$ of shape $D$, we define the {\sl Garnir polynomial}
  $$
  \Delta_T(X)=\prod_{\gamma\in\Gamma(T)}\,\, \det\big\| x_m^{h_T(\ell)}\big\|_{m,\ell\in\gamma}.
  $$
Recall that $h_T(\ell)$ is the height of $\ell$ in $T$, this is to say that it is the first coordinate of the cell in which it
appears in $T$. We have the following proposition.

\begin{prop}
For any tableau $T$ of shape $D'$, we either have
  $${\partial Y}_T\Delta_D(X,Y)=\pm\,\gamma_D\,\Delta_T(X),$$
or ${\partial Y}_T\Delta_D(X,Y)|_{Y=0}=0$. 
\end{prop}

\noindent
Clearly this last possibility can only occur if one of the column of $D$ has a different number
of cells than the corresponding column of $D'$.

\begin{proof}\ We only give an outline of the proof, and restrict ourselves to the case when each column of $D$ has the same number
of cells as the corresponding column of $D'$. In that case,  
$Y_T$ has the same total $Y$-degree as $\Delta_D(X,Y)$, so that ${\partial Y}_T\Delta_D(X,Y)$ is a polynomial
in the $X$ variables only. Since the terms of $\Delta_D(X,Y)$  are
   $${\rm sign(\sigma)}\ x_{\sigma_1}^{r_1}\,y_{\sigma_1}^{c_1}\,x_{\sigma_2}^{r_2}y_{\sigma_2}^{c_2}\,\cdots\,
         x_{\sigma_n}^{r_n}y_{\sigma_n}^{c_n},$$
in view of lemma 2.1, the terms of ${\partial Y}_T\Delta_D(X,Y)$ are forced to be of the form ${\partial Y}_T\,Z_R$
for tableaux $R$ that have the same column set has $T$. On the other hand, since ${\partial
Y}_T\Delta_\mu(X,Y)$ alternates in sign under the action of column fixing permutations, it has to be a multiple of  $\Delta_T(X)$.
Hence, both polynomials having the same degree, we must have the equality stated.
\end{proof}

We now associate to each entry $j$, of a standard tableau $T$, a non negative integer in the following manner.
Let $(r_j,c_j)$ be the position of $j$ in $T$, and let
$k$ be the largest entry of $T$, such that $c_k=c_j+1$ and $k<j$. We set
   $$\alpha(j)=\alpha_T(j):=r_j-r_k.$$
If there is no such $k$, set $\alpha(j):=r_j+1$. For the example given below, the value of $\alpha(k)$ appears in the cell of the
right tableau corresponding to the position of $k$ in the left tableau.
$$\Young{\Box{5}\cr\Box{4}&\Box{8}\cr\Box{3}&\Box{6}\cr\Box{1}&\Box{2}&\Box{7}&\Box{9}&\Box{10}\cr}\qquad
\Young{\Box{3}\cr\Box{2}&\Box{2}\cr\Box{1}&\Box{2}\cr\Box{1}&\Box{1}&\Box{1}&\Box{1}&\Box{1}\cr}
$$      
Clearly, if $T$ is the unique standard tableau corresponding to a column:
$$\Young{\Box{$n$}\cr\hbox{\ \ \vdots}\vskip5pt\cr\Box{2}\cr\Box{1}\cr}
$$ 
then $\alpha(j)=j-1$. 
 
As shown in~\cite{orbit, nantel adriano} the space ${\mathcal L}_{\partial}[\Delta_T(X): T \in
\RS{\mu}]$ 
(the span of all partial derivative of Garnir polynomials for tableaux of shape $\mu$) coincides with the space
$M_\mu^0$. Using this characterization, we will now construct a basis for $M_\mu^0$. But first, let us introduce some further
notation.

For $\mu$ partition of $n$, let $\pi(\mu)$ be the set of partitions of $n-1$ that can be obtained from $\mu$ by removing one of its
corner. For $\nu\in\pi(\mu)$, we denote $\mu/\nu$ the corner by which $\nu$ differs from $\mu$. Let us label $\nu_1,\ldots ,\nu_k$, the partitions in the set $\pi(\mu)$, (for definition, see the proof of theorem 2.3) following
the increasing order of the column number in which the corresponding corners, the $\mu/\nu_i$'s, appear. In other words, if
$(a_i,b_i)$, $1\leq i\leq k$, are the respective coordinates of the corner cells $\mu/\nu_i$, then $b_1<b_2<\ldots < b_k$. Any
standard tableau $T$ of shape $\mu$ is such that $n$ sits in a corner $(a_j,b_j)$ of $\mu$. Moreover, the value of
$\alpha_T(n)$ depends only on the position of this corner (and on the shape $\mu$), since all other entries of $T$ are smaller.
Denoting $\alpha_j$ the value of $\alpha_T(n)$, if $n$ appears in position $(a_j,b_j)$ in $T$, it is clear that
   \begin{equation}\label{alpha j}
       \alpha_j=a_j-a_{j+1}.
   \end{equation}

\begin{thm}\label{basis mu_0} For any partition $\mu$ of n, the set of polynomials
  $${\mathcal B}_\mu:=\big\{\ \partial_ X^{\bf m}\,\Delta_T(X)\ |\  T\in\RS{\mu},\ 
     {\bf m}=(m_1,m_2,\,\ldots\,,m_n) \ {\rm and} \ 0\leq m_i<\alpha_T(i)\ \big\}$$
is a basis of $M_\mu^0$.
\end{thm}

\begin{proof}\ We first show recursively that ${\mathcal B}_\mu$ is independent, assuming that the statement holds for partitions
with at most $n-1$ cells.  As before, for $\nu_j\in\pi(\mu)$, let $(a_j,b_j)$ be the corner $\mu/\nu_j$, with  $b_1<b_2<\ldots <
b_k$, and define
  $${\bf B}_j:=\big\{\ X^{\bf m}\ |\   T\in\RS{\mu},\  \ 0\leq m_i<\alpha_T(i),\ \ 
            T(a_j,b_j)=n\ \big\}.$$
In view of~\ref{alpha j}, for $X^{\bf m}\in {\bf B}_j$,  the dominant monomial of $\partial_ X^{\bf
m}\,\Delta_T(X)$ (in reverse lexicographic order) is of the form
   $$ x_n^{a_j-m_n} X^{\bf p} \qquad {\rm (where}\quad p_n=0{\rm )}, $$
with $0\leq m_n< a_j-a_{j+1} $.
For $k$ fixed with $a_{j+1}<k\leq a_j$,  our induction hypothesis gives that the set
  $${\mathcal B}_{j,k}:= \big\{  \partial_ X^{\bf m}\,\Delta_T(X)\ |\   T\in\RS{\mu},\  \ 0\leq m_i<\alpha_T(i),\ \ 
            T(a_j,b_j)=n,\ m_n=a_j-k\ \big\}\,,$$
is independent, since (in reverse lexicographic order) we have the following expansion of $\partial_ X^{\bf m}\,\Delta_T(X)$ 
  $$\partial_ X^{\bf m}\,\Delta_T(X)=x_n^k \partial_ X^{\bf p}\,\Delta_{T'}(X) + \underbrace{\ldots}_{\rm lower\ terms}$$
where $T'$ is the restriction of $T$ to $\nu_j$. Clearly, the sets ${\mathcal B}_{j,k}$ are mutually independent, so
   \begin{equation}\label{base mu}
      {\mathcal B}_\mu=\biguplus_{j,k}{\mathcal B}_{j,k}
   \end{equation}
is independent.

We will now show that the number of elements of ${\mathcal B}_\mu$ is 
   \begin{equation}\label{dimension mu}
      | \, {\mathcal B}_\mu\,|= {n!\over \mu!}\,,
   \end{equation}
using a recursive argument, assuming that the statement holds for $\nu_j\in\pi(\mu)$. By induction, it is clear that
    $$ | \, {\mathcal B}_{j,k}\,|= {(n-1)!\over \nu_j!}$$
so that (in view of~\ref{base mu})
  $$| {\bf B}_\mu | = \sum_{i=1}^k \alpha_j \, {(n-1)!\over \nu_j!}\,.$$
The result follows from the easy observation that
  $$n=\sum_{j=1}^k  \alpha_j\,(a_j+1), $$
since 
   $$a_j+1={\mu!\over \nu_j!}$$ 
is the length of the row of $\mu$ in which sits the corner $(a_j,b_j)$.
A "geometric" argument, that can be found in \cite{garsia procesi}, shows that the dimension of  $M_\mu^0$ is  at most
$\frac{n!}{\mu!}$. Thus ${\mathcal B}_\mu$ is a basis.
\end{proof}

It is shown in \cite{garsia procesi} that the graded Frobenius character of the $M_\mu^0$'s are none other than the
Hall-Littlewood symmetric functions.

\section{A Basis for $M_{\mu/ij}^0$}

The central result of this paper is the following description of a basis for 
$M_{\mu/i,j}^0$, with $\mu$ any fixed partition of $n+1$, and $(i,j)$ any given cell of $\mu$. The proof that it is a generator set
is postponed until the next section. 

For a standard tableau $T$, let ${\bf B}_T$ simply denote the set
   $${\bf B}_T:=\big\{\ X^{\bf m}\ |\  0\leq m_s\leq \alpha_T(s)\ \big\}.$$
For $\nu_\ell$ a partition of $n$ obtained from $\mu$ by removing the corner cell $(a_\ell,b_\ell)$, the basis of
$M_\ell:=M_{\nu_\ell}$, described in the previous section, is
   $${\mathcal B}_\ell=\big\{\partial_ X^{\bf m}\Delta_T(X)\ |\ T\in\RS{\nu_\ell},\   X^{\bf m}\in {\bf B}_T\ \big\}.$$ 
If $T$ is a standard tableau of shape $\nu_\ell$, and $0\leq u\leq a_\ell$ an integer, we denote
$T\!\uparrow_{uv}$ the tableau of shape $\mu/uv$ (with  $v=b_\ell$), such that
   $$T\!\uparrow_{uv}(r,c)=\left\{
\begin{array}{cl}
T(r,c)\hfill & \ \hbox{if} \ c\not=v\ \hbox{or}\  r< u\,, \\
                 \\
T(r-1,c)  & \ \hbox{if} \ c=v\ \hbox{and}\ r> u \,.
\end{array}
\right.$$
Since $(u,v)$ is not in $\mu/uv$, $T\!\uparrow_{uv}$ need not be defined at $(u,v)$. In other words, the tableau $T\!\uparrow_{uv}$
is obtained from $T$ by ``sliding'' upward by $1$ the cells in column $v$ that are on or above row $u$. For $u$ and $v$ as above,
we set
   \begin{equation}\label{partie base}
      {\mathcal A}_{uv}:=\big\{ \partial_ X^{\bf m} \Delta_{T\uparrow_{u,v}}(X)\ |\ 
         X^{\bf m}\in{\bf B}_T,\ T\in \RS{\nu_\ell}    \big\}.
   \end{equation}
Observe that ${\mathcal A}_{uv}$ implicitly depends on the choice of corner $(a_\ell,b_\ell)$ of $\mu$, since
$v$ is  equal to  $b_\ell$. Moreover,   ${\mathcal A}_{a_\ell,b_\ell}$ is the basis of $M_{\nu_\ell}^0$
described in theorem~\ref{basis mu_0}, thus ${\mathcal A}_{a_\ell,b_\ell}$ is independent.

Let $(a_1,b_1)$, $\ldots$ ,$(a_m,b_m)$ (with $b_1<\cdots<b_m $) be the set of corners of $\mu$ that are in the ``shadow'' of
$(i,j)$. This is to say that  $i\leq a_\ell$ and $j\leq b_\ell$, for all $1\leq \ell\leq m$. Once again, we denote 
$\alpha_\ell$  the value of
$\alpha_T(n+1)$ for any standard tableau of shape $\mu$ having $n+1$ in position
$(a_\ell,b_\ell)$. 
Defining 
   \begin{equation}\label{base muij}
       {\mathcal B}_{\mu/ij}:=\bigcup_{\ell=1}^m \bigcup_{u=i}^{\min(i+\alpha_\ell-1,a_\ell)} {\mathcal A}_{u,b_\ell}\,,
   \end{equation}
we have 

\begin{thm}\label{rec basis}
  For $\mu$ a partition of $n+1$ and $(i,j)\in\mu$, ${\mathcal B}_{\mu/ij}$ is a basis of $M_{\mu/ij}^0$.
\end{thm}

\begin{proof}\ In the remainder of this section, we will prove that~\ref{base muij} is an independent set,
using a downward recursive argument, and that the number of elements of  ${\mathcal B}_{\mu/ij}$ is 
   \begin{equation}\label{dimension muij}
       d_{\mu/ij}:={\frac {n!} {\mu!}} \sum_{\twolines {i'> i} {\mu_{i'}> j}} \mu_{i'}\,.
   \end{equation}
To complete the proof of the
theorem, we will show in section 4 that the dimension of
$M_{\mu/ij}^0$ is at most equal to  $d_{\mu/ij}$, so that $M_{\mu/ij}^0$ has to coincide with the span of ${\mathcal B}_{\mu/ij}$.
\end{proof}

The cardinality of ${\mathcal A}_{u,b_\ell}$ is clearly
  $${\frac {n!} {\mu!}}\,\mu_{a_\ell+1-u+i} ={\frac {n!} {\mu!}}\,\mu_{a_\ell+1} $$
since this is the same as ${ {n!}/{\nu_\ell!}}$. If~\ref{base muij} is a disjoint union (which would follow from it
being independent) and since every $i'$ indexing the summation in~\ref{dimension muij} occur
exactly once in $\{a_\ell+1-u+i\ |\ 1\le\ell\le m,\ i\le u\le
\min(i+\alpha_\ell-1,a_\ell)\}$, then we must have
   \begin{equation}\label{egalite}
       |{\mathcal B}_{\mu/ij}|=d_{\mu/ij}\,.
   \end{equation}
Let
   $$D_X:=\partial_{x_1}+\partial_{x_2}+\ldots + \partial_{x_n}\,.$$
We recall the following special case of proposition (I.2) of \cite{Berg et al}.

\begin{prop} If $(i+1,j)$ is in $\mu$, then 
  $$D_X\,\Delta_{\mu/i,j}(X,Y)=_{\rm cte} \Delta_{\mu/i+1,j}(X,Y).$$
Otherwise $D_X\,\Delta_{\mu/ij}(X,Y)=0$, and this corresponds to the case where $(i,j)$ is on the top border of $\mu$. 
The symbol ``$=_{\rm cte}$'' stands for equality up to a non zero constant.
\end{prop}

It is easy to adapt the proof of this fact to show that, for $T$  a standard tableau of shape
$\nu_\ell$ (with $\mu/\nu_\ell=(a_\ell,b_\ell)$) and
$0\leq u\leq a_\ell$, we have
   \begin{equation}\label{DX tableau}
       D_X\, \Delta_{T\uparrow_{u,v}}(X)=_{\rm cte}\left\{
\begin{array}{cl}
\Delta_{T\uparrow_{u+1,v}}(X)\hfill & \ \hbox{if} \ u<a_\ell\,, \\
                 \\
0  & \ \hbox{if} \ u=a_\ell \,,
\end{array}
\right.
   \end{equation}
where  $v=b_\ell$. It follows from definition~\ref{partie base} that

\begin{lem}\label{Dx} Using the same convention as above, we have
   \begin{equation}\label{DX base partie}
       D_X\,{\mathcal A}_{u,b_\ell}=_{\rm cte}\left\{
\begin{array}{cl}
{\mathcal A}_{u+1,b_\ell}\hfill & \ \hbox{if} \ u<a_\ell\,, \\
                 \\
\{0\}  & \ \hbox{if} \ u=a_\ell \,.
\end{array}
\right.
   \end{equation}
\end{lem}
Since these two sets have the same cardinality we deduce, from the linear independence of 
${\mathcal A}_{a_\ell,b_\ell}$, that each ${\mathcal A}_{u,b_\ell}$ is independent. Applying $D_X$ in
definition~\ref{base muij} we readily check that
   \begin{equation}\label{DX base}
      {\mathcal B}_{\mu/i+1,j}=D_X\,{\mathcal B}_{\mu/i,j}\,.
   \end{equation}
But we know that $D_X\,{\mathcal A}_{a_m,b_m}=\{0\}$, and it is clear that ${\mathcal A}_{a_m,b_m}$ is a subset of
${\mathcal B}_{\mu/ij}$. A dimension count, together with the recursive assumption, forces ${\mathcal B}_{\mu/i+1,j}$
to be independent, since
   $$d_{\mu/i+1,j}+|{\mathcal A}_{a_m,b_m}| ={\frac {n!} {\mu!}} \sum_{\twolines {i'> i+1} {\mu_{i'}> j}}
\mu_{i'}+{\frac {n!} {\mu!}}\,\mu_{i+1}=d_{\mu/ij}\,.
$$
This ends the proof of theorem~\ref{rec basis}.
\section{Upper bound for the dimension of $M_{\mu/ij}^0$ }

We now give an upper bound for the dimension of $M_{\mu/ij}^0$.
Given a polynomial $P(X,Y)$, we denote by $P(\partial)$ the operator obtained
 from $P$ by replacing all
 the variables $x_i$ and $y_j$ by $\partial_{x_i}$ and $\partial_{y_j}$,
respectively.
If $P(X,Y)=M$ is a single monomial $M$ we will write $P(\partial)=\partial_M$.

\begin{thm}\label{th:bound}
  For $\mu$ a partition of $n+1$,
    $$\dim M_{\mu/i,j}^0 \le {\frac {n!} {\mu!}} \sum_{\twolines {i'> i}
       {\mu_{i'}> j}} \mu_{i'}.$$
\end{thm}

\begin{proof}\ In \cite{Berg et al}, the bigraded $\S_n$-modules
$M_{\mu/i,j}$ and
${\mathcal
L}_{\partial}[\partial_{x_{n+1}}^i\partial_{y_{n+1}}^j\Delta_{\mu}(X;Y)]$
are shown to be equivalent, hence their $Y$-free components
  $$M_{\mu/i,j}^0 \qquad {\rm and}\qquad {\mathcal
L}_{\partial}[\partial_{x_{n+1}}^i\partial_{y_{n+1}}^j\Delta_{\mu}(X;Y)]^0\,,$$
are equivalent. For any injective tableau $T$ of shape $\mu$, with
$c_{n+1}\ge j$, we have
$Y_T=y_{n+1}^j Y^{\bf a}$ and
  $$\pm\gamma_\mu\,\Delta_T(X)=\partial Y_T\Delta_{\mu}(X,Y).$$
If $r_{n+1}< i$ then $\partial_{x_{n+1}}^i \Delta_T(X)=0$,
so that ${\mathcal
L}_{\partial}[\partial_{x_{n+1}}^i\partial_{y_{n+1}}^j\Delta_{\mu}(X;Y)]^0$
is equal
to
   \begin{equation}\label{Eq:n}  {\mathcal
L}_{\partial}\big[\partial_{x_{n+1}}^i \Delta_T(X)\ \big|\
         T\colon\mu\to\{1,2,\ldots,n+1\},\
          c_{n+1}\ge j,\ r_{n+1}\ge i\big]\,.
   \end{equation}

For $\zeta_1,\zeta_2,\ldots,\zeta_{\ell(\mu)}$  and
$\omega_1,\omega_2,\ldots,\omega_{\ell(\mu')}$  two families of pairwise
distinct scalars, we construct a set of points
$[\rho_{\mu}]_{i,j}$  in $\C^{2(n+1)}$ as follow.
For every injective tableau $T$ of shape $\mu$, we
define the point
$\rho_{_T}=(\zeta_{r_1},\zeta_{r_2},\ldots,\zeta_{r_{n+1}},\omega_{c_1},
            \omega_{c_2},\ldots,\omega_{c_{n+1}})$
in $\C^{2(n+1)}$, and set
  $$[\rho_{\mu}]_{i,j}=\big\{\rho_{_T}\ \big|\
T\colon\mu\to\{1,2,\ldots,n+1\},\ c_{n+1}\ge j,\
r_{n+1}\ge
      i\big\}.$$
That is $\rho_{_T}\in[\rho_{\mu}]_{i,j}$ when $n+1$ lies in the shadow of
$(i,j)$ in $T$. Note that
$\rho_{_T}\in[\rho_{\mu}]=\rho_{_T}\in[\rho_{\mu}]_{0,0}$ contain $n!$
points in correspondance with
every injective tableau of shape $\mu$.

We denote by $[\rho_{\mu}]_{i,j}^0=\pi([\rho_{\mu}]_{i,j})$, where $\pi$ is
the projection on
$\C^{n+1}$ that keeps only the first $n+1$ entries.
We see that the set of tableaux with
$n+1$ entries strictly increasing in rows and where $n+1$ lies in a row
$i'$ such that $i'=r_{n+1}+1> i$ and $\mu _{i'}=c_{n+1}+1> j$ give
all the points of $[\rho_{\mu}]_{i,j}^0$ exactly once.
One then easily verifies that the cardinality of
$[\rho_{\mu}]^0_{i,j}$ is precisely
  $$d_{\mu/ij}={\frac {n!} {\mu!}}
   \sum_{\twolines {i'> i} {\mu_{i'}> j}} \mu_{i'}.$$
Following \cite[Section 4]{Berg et al} we associate to
this set $J_{[\rho_{\mu}]^0_{i,j}}$, its annulator ideal and define
$H_{[\rho_{\mu}]^0_{i,j}}=(J_{[\rho_{\mu}]^0_{i,j}})^{\perp}$. The
dimension of $H_{[\rho_{\mu}]^0_{i,j}}$ is then $d_{\mu/ij}$ as well.

Given a polynomial $P$, let $h(P)$ denotes its  homogeneous
component of highest degree. For any polynomial $P$ in
$J_{[\rho_{\mu}]^0_{i,j}}$, let
  $$
    Q(X,Y)=P(X)\prod_{i'=1}^i(x_{n+1}-\zeta_{i'})
    \prod_{j'=1}^j(y_{n+1}-\omega_{i'})
  $$
For any $\rho_{_T}\in [\rho_{\mu}]$, the two products in the definition of
$Q$ vanish at $\rho_{_T}$
unless
$n+1$ lies in the shadow of $(i,j)$ in $T$. But if this is the case then
$P(X)$ vanishes at
$\pi(\rho_{_T})$. This shows that $Q(X,Y)$ is in
$J_{[\rho_{\mu}]}$, the annulator ideal of
$[\rho_{\mu}]_{0,0}=[\rho_{\mu}]$.  Hence
$h(Q)=h(P)x_{n+1}^iy_{n+1}^j$ is in ${\rm gr}(J_{[\rho_{\mu}]})$, its graded
version and
$h(Q)(\partial)\Delta_\mu(X,Y)=0$.  For any injective tableau
$T$ of shape
$\mu$ such that $c_{n+1}\ge j$ and $r_{n+1}\ge i$ we have $Y_T=y_{n+1}^j
Y^{\bf a}$ and
  \begin{eqnarray*}
   h(P)(\partial)\,\partial_{x_{n+1}}^i\Delta_T(X)
   &=& \pm \gamma_\mu^{-1}
\,h(P)(\partial)\,\partial_{x_{n+1}}^i\,\partial Y_T\Delta_{\mu}(X,Y)\\
   &=& \pm \gamma_\mu^{-1}
\,h(P)(\partial)\,\partial_{x_{n+1}}^i\,\partial_{y_{n+1}}^j
        \partial_Y^{\bf a}\Delta_{\mu}(X,Y)\\
   &=& \pm \gamma_\mu^{-1}\partial_Y^{\bf a}
\,h(Q)(\partial)
        \Delta_{\mu}(X,Y)\  = \ 0 \\
\end{eqnarray*}
Thus $h(P)$ is in
$I_{\partial_{x_{n+1}}^i\,\Delta_T(X)}$. We obtain this way that
gr$J_{[\rho_{\mu}]^0_{i,j}}$ is a
subset of
$I_{\partial_{x_{n+1}}^i\,\Delta_T(X)}$ for any
$T$ with the prescribed conditions. The space in~\ref{Eq:n} is thus
contained in
$H_{[\rho_{\mu}]^0_{i,j}}$, which proves the theorem.
\end{proof}

\section{Four term recurrence}

Specializing conjecture~\ref{conj:rec} to its
$Y$-free component, corresponds to setting $q=0$ in the four term recurrence~\ref{quatres termes}. We now show
that this specialization of conjecture~\ref{conj:rec} holds, by giving an explicit interpretation of the
resulting recurrence in term of the basis we have constructed for $M_{\mu/ij}^0$. 
 
\begin{thm}\label{4termrec}
If $H_{\mu/i,j}^0$ denotes the graded Frobenius characteristic of $M_{\mu/ij}^0$ then:
\begin{itemize}
\item if $a=0$ and $\ell>0$,\quad $\displaystyle H_{\mu/i,j}^0=\frac {t^{\ell+1}-1}
{t^{\ell}-1}\,H_{\mu/i,j+1}^0$;
\medskip
\item if $a>0$,\quad
$\displaystyle H_{\mu/i,j}^0=H_{\mu/i,j+1}^0+t\,H_{\mu/i+1,j}^0-t\,H_{\mu/i+1,j+1}^0$.
\medskip
\item if $a=0$ and $\ell=0$,\quad $H_{\mu/ij}^0$ is the graded Frobenius characteristic of $M_{\nu}^0$, where $\nu$ is the
partition $\mu/ij$.
\end{itemize}
Here (as before) $\ell$ and $a$  give the number of cells that are respectively 
north and east of
$(i,j)$ in  $\mu$. If any of the cells $(i+1,j)$, $(i,j+1)$ or $(i+1,j+1)$ falls out of $\mu$, then the corresponding term is
considered to be 0.
\end{thm}

\begin{proof}\ Each of these assertions can be shown using the basis we have constructed. The third one is just a direct
observation. The first one corresponds to a case for which there is just one corner $(a_m,b_m)$ in the shadow of $(i,j)$, with
$b_m=j$, and then~\ref{base muij} can be written as
   $${\mathcal B}_{\mu/ij}:= \bigcup_{u=0}^{\ell} {\mathcal A}_{i+u,j}\,,$$
Since, as long as $(k+1,j)$ is in $\mu$, $D_X$ is an isomorphism of representations between the homogeneous $\S_n$-modules
${\mathcal A}_{k,j}$ and ${\mathcal A}_{k+1,j}$ that lowers the degree by 1, we must have
  $${\mathcal F}_t({\mathcal A}_{k,j})=t\,{\mathcal F}_t({\mathcal A}_{k+1,j}) $$
where ${\mathcal F}_t$ stands for the graded Frobenius characteristic. We deduce that, in the first case,
   $${\mathcal B}_{\mu/ij}=(1+t+\ldots +t^\ell)\, H_\nu^0$$
with $\mu/\nu=(a_m,b_m)$. This is clearly equivalent to the statement of the first case.

For the second case there are a few subcases, all similarly dealt with, the most interesting one being when $j=b_1$ and $m>1$ for
which the basis can clearly be broken down as
  $${\mathcal B}_{\mu/ij}={\mathcal B}_{\mu/i,j+1} \uplus \bigcup_{u=i}^{ i+\alpha_1-1 } {\mathcal A}_{u,b_1}$$
and we only need to show that the graded Frobenius characteristic of the linear span of 
  $$\bigcup_{u=i}^{ i+\alpha_1-1 } {\mathcal A}_{u,b_1}$$
is given by 
   \begin{equation}\label{noyau}
      t\,(H_{\mu/i+1,j}^0-H_{\mu/i+1,j+1}^0).
   \end{equation}
Now we clearly have ${\mathcal B}_{\mu/i+1,j+1}\subset {\mathcal B}_{\mu/i+1,j}$, with $D_X\,\bigcup_{u=i}^{
i+\alpha_1-1 } {\mathcal A}_{u,b_1}$ being the complement of ${\mathcal B}_{\mu/i+1,j+1}$ in ${\mathcal
B}_{\mu/i+1,j}$. Under the hypothesis of this subcase, the graded Frobenius characteristic of the span of
$\bigcup_{u=i}^{ i+\alpha_1-1 } {\mathcal A}_{u,b_1}$ is thus given by~\ref{noyau}. All other subcases are simple to
show.   
\end{proof}

\section{Remarks }

\begin{Remark} {\em In \cite{nantel adriano} (proposition 2.2), N. Bergeron and Garsia show that the spaces $M_\mu^0$ are nested
into each other according to their partition indexing. That is
  $$\mu\preceq\lambda\qquad\implies\qquad M_\mu^0\subseteq M_\lambda^0 $$
where ``$\preceq$'' denotes the dominance order. Moreover they show that
  $$ M_\mu^0\cap M_\lambda^0 =M_{\mu\wedge\lambda}^0.$$
Using our basis, it is easy to show that both these results extend to the situation studied in this paper. Namely,}
\end{Remark}

\begin{prop}
For two partition $\mu$ and $\lambda$ of $n+1$, we have
  $$\mu\preceq\lambda\qquad\implies\qquad M_{\mu/ij}^0\subseteq M_{\lambda/ij}^0 $$
and
  $$ M_{\mu/ij}^0\cap M_{\lambda/ij}^0 =M_{\mu\wedge\lambda/ij}^0,$$
whenever $(i,j)$ appears in both $\mu$ and $\lambda$.
\end{prop}

\begin{Remark} {\em For $\mu$ a partition of $n$ (denoted $\mu\vdash n$), Macdonald has given an explicit description of the
coefficients appearing in the Pieri formula for the $H_\mu$:
  $$h_k^\perp\, H_\mu(X;q,t)=\sum_{\nu\vdash n-k\atop \nu\subseteq \rho}c_{\mu\nu}^k(q,t) H_\nu(X;q,t)$$
where $h_k^\perp$ is the operator dual to multiplication by $h_k$ (complete homogeneous) with respect to the usual scalar product
on symmetric functions for which the Schur functions are orthonormal. These coefficients $c_{\mu\nu}^k(q,t)$ are rational functions
in $q$ and $t$. Now, let $\rho$ be the partition of $m$ corresponding to the shadow of $(i,j)$ in $\mu$, with $m$ equal to the
number of cells in this shadow. F. Bergeron has conjectured in \cite{bergeron} that the following symmetric function
   \begin{equation}\label{formule pieri}
       \sum_{\nu\vdash m-k\atop \nu\subseteq \rho}c_{\mu\nu}^k(q,t) H_{\mu-\rho+\nu}(X;q,t),
   \end{equation}
where $\mu-\rho+\nu$ stands for the partition obtained from $\mu$ by replacing $\rho$ (the shadow of $(i,j)$) by $\nu$,
is the bigraded Frobenius characteristic of the module $M_{\mu/ij}^k$ obtained as the union of all modules $M_D$, for
$D$ ranging in the set of diagrams obtained from $\mu$ by removing $k$ cells in the shadow of $(i,j)$. This would imply
that the dimension of 
$M_{\mu/ij}^k$ be equal to ${m\choose k}\,(n-k)!$. J.-C. Aval, in \cite{aval certain}, has shown that this value is an upper
bound, and has generalized the construction of this paper to obtain an explicit basis for the $Y$-free component of
$M_{\mu/ij}^k$. One can show that the graded Frobenius characteristic of the resulting space is the symmetric function
obtained by taking the limit as $q\rightarrow 0$ of~\ref{formule pieri}.}
\end{Remark}

\begin{Remark} {\em 
One can explicitly characterize the defining ideal of the space $M_{\mu/ij}^0$. This will be the subject of a forthcoming paper~\cite{aval tani}.
}\end{Remark}

\end{document}